\documentclass[12pt,a4paper]{article}
\usepackage{authblk,amsmath,amsfonts,amssymb,amsthm,graphicx,caption,subcaption}
\usepackage[T2A]{fontenc}
\usepackage[utf8]{inputenc}
\usepackage[english]{babel}
\usepackage[top=2.0cm,bottom=2.0cm,left=2.0cm,right=2.0cm]{geometry}
\usepackage{multirow}

\usepackage[noend]{algorithmic}
\usepackage{algorithm}
\usepackage{booktabs}

\providecommand{\keywords}[1]{\textbf{\textit{Keywords:}} #1}
\newcommand\fundnote[1]{\textit{#1}}

\newcommand{\figref}[1]{\figurename~\ref{#1}}

\numberwithin{equation}{subsection}

\begin{document}

\begin{center}
\Huge{Acquiring elastic properties of thin composite structure from vibrational testing data}\\[1cm]
\Large{V. Aksenov, A. Vasyukov, K. Beklemysheva}\\[0cm]
\large{MIPT}\\[0cm]
\large{141700, Moscow region, Dolgoprudny, Institutsky lane, 9}\\[0cm]
\large{E-mail: aksenov.vv@phystech.edu}\\[0cm]
\end{center}

\abstract{%
    The problem of acquiring elastic properties of a composite material from the data of the vibrational testing is considered. 
    The specimen is considered to abide by the linear elasticity laws and subject to viscous-like damping.
    The BVP for transverse movement of such a specimen under harmonic load is formulated and solved with finite-element method.
    The problem of acquiring the elastic parameters is then formulated as a nonlinear least-square optimization problem.
    The usage of the automatic differentiation technique for stable and efficient computation of the gradient and hessian allows to use well-studied first and second order optimization methods, namely Newton's and BFGS.
    The results of the numerical experiments on simulated data are analyzed in order to provide insights for the experiment planning.
}


\fundnote{This study is supported by RSF grant 19-71-00147}

\bigskip
\section{Introduction}

Due to their advanced properties, composite materials are widely used in modern-day engineering projects, and new materials, possibly tailor-made for a certain task, appear regularly.
However, the experimental status of such materials may lead to lack of reliable data on their elastic properties, crucial for numerical simulations.
The aforementioned properties may also change in the process of production of a certain component.
Another feature of composite materials is that under certain extreme conditions they can develop different types of internal damage, invisible from the outside, but detrimental for the operational characteristics. 
Thus, the need for nondestructive  experimental  methods of studying the properties of the particular composite component and reliable numerical algorithms for the experimental data processing arises.

This work focuses on composite plates, i.e. specimen for which a spatial direction (denoted as $z$ direction) exists, in which the linear size of the specimen is considerably smaller than in two other directions.
One of the known applications of composite plates is the production of electrode for ion thrusters.
The research in this area is conducted by numerous groups around the world, for example \cite{MadeevThruster, nishiyama2005research}.

\subsection{Previous theoretical and experimental studies}
    In \cite{tuan2015exploring} the similar experimental stand is considered. 
    The mechanical impedance of a triangular aluminum plate, driven by a shaker at different frequencies, is measured.
    Experiments are conducted both in the atmosphere and in the vacuum chamber, and the difference between the frequency response is minor, indicating that the main contribution to the energy dissipation is due to the damping properties of the material itself rather than the interaction with the external medium. The internal damping processes, such as thermo- or viscoelasticity, should be carefully studied.
    Another important finding is the discrepancy between the natural frequencies of the specimen, obtained from the corresponding eigenvalue problem, and the extrema of the empiric frequency response function.
    
    In \cite{treviso2015damping} the different models for damping in composite materials are overviewed. 
    The most prominent model, used in the analysis of composite plates, works with the Laplace transform over on the original dynamic linear elastic equations and allows the elastic moduli to take complex values, i.e. $\hat{D}_{ij} = D_{ij}(1 + i\beta_{ij})$.
    Here $D_{ij}$ is the static elastic modulus from the generalized Hooke's law, and $\beta_{ij}$ is the corresponding loss factor.
    These coefficients do, in principle, depend on the frequency \cite{pritz1998frequency}.
    The powerful tool for studying the dynamic moduli is Dynamic Mechanical Analysis \cite{menard2020dynamic, ayyagari2019enhancing}.
    This experimental method, however, requires small samples of the material, possessing certain shapes, making it useless when we assume that the elastic properties are not consistent between specimens.
    
    An extensive amount of publications consider the study of the elastic properties by solving the problems, posed for a model of a thin plate, in which the displacement field abides by the kinematic hypotheses of either Kirchhoff-Love or Reissner-Mindlin. 
    The governing equations for these models can be found in \cite{reddy2003mechanics}.
    In several of such publications, the damping is either not considered at all \cite{pagnotta2008elastic,lee2006identification, li2019vibration, barkanov2015effectiveness}, or is considered to be frequency-independent in the frequency range in consideration \cite{qian1997vibration,schwaar2012modal}. This, as noted by the authors in \cite{matter2009numerical}, may lead to the discrepancy between experimental and numerical modal frequencies and decay factors for higher frequencies in composites with viscous polymer constituents. The usage of frequency-dependent dynamic moduli in \cite{ruzek2013experimental} provided better explanation of the dynamic behaviour of a specimen made of two steel plates connected by a polymer layer.
    In \cite{chaigne2001time} a framework is suggested, which allows to introduce into the model the phenomena of viscous, thermoelastic and viscoelastic damping and the damping due to sound radiation. 
    The asymptotic models for the complex moduli, used by the authors, allowed for good reconstruction of the experimental natural frequencies and decay factor of specimens made of aluminum, glass, carbon fiber and wood.
    Most of the previous studies (with exception of \cite{sankar2014system}, which uses time-domain data) use modal frequencies and decay factors and mode shapes as the data for reconstruction. 
    These, in turn, are obtained via a pipeline of the scanning vibrometer and a modal post-processing software. 
    The authors also mostly consider rectangular plates made of orthotropic, but not of more general monoclinic material.

    
\section{Experimental setup and physical assumptions} 
    \subsection{Description of experimental stand}
    The composite plate specimen is attached with bolts to a special stand which can vibrate at a given frequency, controlled by the experimenter.
    Two accelerometers, one attached to the stand, and another~---~to the surface of the specimen, are synchronized and measure the dependence of the acceleration at each point as a function of time. 
    Suggesting that both the movement of the stand and the response in the test point are harmonic vibrations at the driving frequency, amplitude ratio and phase shift between these vibrations can be measured from these data.
    Amplitude ratio and phase shift can be alternatively expressed as a complex number, and its dependence on the driving frequency is called the \emph{amplitude-frequency characteristic} in the following sections.
    In the sections below the authors explain how to reproduce the results of such experiment for a specimen with given geometry and physical properties  in a numerical simulation and explore the possibility of reconstructing the material properties from the AFC.
    \subsection{Governing equations}
    The material occupies the volume $\overline{\Omega} \times [-\frac12 h, \frac12 h] \subset \mathbb{R}^3$, where $h$ is the plate thickness,
     and the open connected set $\Omega$ is the midplane of the plate. 
    The boundary of the plate comprises of two parts: $\partial \Omega = \Gamma_c \sqcup \Gamma_f$, 
     with $\Gamma_c$ being the clamped edge, and $\Gamma_f$ the free edge.

    It is believed that the displacement field in the volume of the thin plate can be defined as a function of the mid-plane 
        displacement $w^0$ via the \emph{Kirchhoff-Love kinematic hypothesis}
    \begin{align*}\label{eq:KirchhoffLove}
        u(x, y, z) &= -z\frac{\partial w^0(x, y)}{\partial x} \\
        v(x, y, z) &= -z\frac{\partial w^0(x, y)}{\partial y} \\
        w(x, y, z) &= w^0(x, y)
    \end{align*}
    Plugging this form of the displacement field into the laws of motion and integrating through thickness leads to the equation (adapted from \cite{bonaldi2018mathematical})
    \begin{equation*}
        2e\rho \ddot{w}^0 - \frac{2e^3}{3}\rho\Delta\ddot{w}^0 - \operatorname{div}\operatorname{div}\mathbf{M} = f
    \end{equation*}
    Here $e = \frac12h$, $\rho$ is the material density, $f$ is the external loading force, and the tensor of moments $\mathbf{M}$ is defined as
    \begin{equation*}
        \begin{pmatrix}
            M_{xx} & M_{xy} \\
            M_{xy} & M_{yy}
        \end{pmatrix} =
        \int\limits_{z = -\frac12h}^{z = \frac12h}z^2
        \begin{pmatrix}
            \sigma_{xx} & \sigma_{xy} \\
            \sigma_{xy} & \sigma_{yy}
        \end{pmatrix} 
    \end{equation*}
    In the present work, the material is considered to abide by the laws of linear elasticiy, with elastic parameters constant through the thickness of the plate.
    This can be expressed as
    \begin{equation}
        \begin{pmatrix}
            M_{xx} \\
            M_{yy} \\
            M_{xy}
        \end{pmatrix} = 
        \begin{pmatrix}
            D_{11} & D_{12} & D_{16} \\
            D_{12} & D_{22} & D_{26} \\
            D_{16} & D_{26} & D_{66} 
        \end{pmatrix}
        \begin{pmatrix}
            \frac{\partial^2}{\partial x^2}w^0 \\
            \frac{\partial^2}{\partial y^2}w^0 \\
            \frac{\partial^2}{\partial x \partial y}w^0 
        \end{pmatrix}=\frac{2e^3}{3}
        \begin{pmatrix}
            C_{11} & C_{12} & C_{16} \\
            C_{12} & C_{22} & C_{26} \\
            C_{16} & C_{26} & C_{66} 
        \end{pmatrix}
        \begin{pmatrix}
            \frac{\partial^2}{\partial x^2}w^0 \\
            \frac{\partial^2}{\partial y^2}w^0 \\
            \frac{\partial^2}{\partial x \partial y}w^0 
        \end{pmatrix}
    \end{equation}
    with $C_{ij}$ being the constants from the generalized Hooke's law. 
    For the isotropic material, 
    \begin{equation}\label{eq:isotropicHooke}
        \begin{pmatrix}
            D_{11} & D_{12} & D_{16} \\
            D_{12} & D_{22} & D_{26} \\
            D_{16} & D_{26} & D_{66} 
        \end{pmatrix} =
        \frac{2Ee^3}{3(1-\nu^2)}
        \begin{pmatrix}
            1 & \nu & 0 \\
            \nu & 1 & 0 \\
            0 & 0 & 1 - \nu 
        \end{pmatrix} =
        D
        \begin{pmatrix}
            1 & \nu & 0 \\
            \nu & 1 & 0 \\
            0 & 0 & 1 - \nu 
        \end{pmatrix}
    \end{equation}
    where $D$ is the \emph{flexural rigidity} of the plate.
    
    Taking into account the above and dividing the equation by $2e$, we arrive at the following transient problem, completed by the boundary conditions.
    \begin{align}\label{eq:transient_start}
        &\rho \ddot{w}^0 - \frac13 \rho e^2\Delta\ddot{w}^0 - \frac{1}{2e}\operatorname{div}\operatorname{div}\mathcal{D}\nabla\nabla w^0 = \frac{1}{2e}f \\
        &\begin{cases}
            w^0 = g \\
            \frac{\partial w^0}{\partial n} = 0
        \end{cases} on\ \Gamma_c \\
        &\begin{cases}
            \frac23 e^3 \rho \frac{\partial}{\partial n}\ddot{w}^0 + \operatorname{div}\mathbf{M \cdot n} + \frac{\partial}{\partial \tau}(\mathbf{Mn \cdot \tau}) = 0 \\
            \mathbf{Mn \cdot \tau} = 0
        \end{cases} on\ \Gamma_f 
        \label{eq:transient_fin}
    \end{align}
    Here $\mathcal{D}$ consists of $D_{ij}$ but they are arranged as a symmetric rank $4$ tensor (for correct notation in this equation).
\section{Mathematical model} 
    \subsection{BVP for AFC}
        In the experimental setup described above, the clamped part of the boundary vibrates with known amplitude and frequency:
        \begin{equation*}
            g(x, y, t) = g_\omega(x, y) \cdot e^{i\omega t}
        \end{equation*}
        Due to the BVP \eqref{eq:transient_start}-\eqref{eq:transient_fin} being linear, it is possible to represent the solution as
        \begin{equation*}
            w^0 = w^0_{part} + \sum\limits_{i = 0}^{\infty}C_iw^0_i
        \end{equation*} 
        with the second term being the linear combination of solutions of the BVP with zero boundary condition that satisfies some initial conditions.
        Due to the presence of damping it is reasonable to assume that these terms fade exponentially, thus, if the measurement is done after a long enough time after the initiation of the vibrations, the response will be roughly the same for all possible initial conditions.
        
        The partial solution is found via separation of variables:
        \begin{equation*}
            w_{part} = u(x, y, \omega)\cdot e^{i\omega t}
        \end{equation*}
        Here $u$ can be viewed as the complex amplitude of the driven vibrations, with $|u|$ representing the amplitude and $\arg{u}$~---~the phase shift of the driven vibrations in relation to the phase of the driving vibration.

        On this stage we also introduce the damping in the equations by substituting each of the elastic parameters with a complex value
        \begin{equation}\label{eq:damping}
            \hat{D}_{\alpha} = D_{\alpha}(1 + i\beta_{\alpha})
        \end{equation}
        where $D_{\alpha}$ is called \emph{storage modulus} and $\beta_{\alpha}$ is called \emph{loss factor} ($\alpha \in \{11, 12, 16, 22, 26, 66 \}$).

        We thus arrive at the BVP for finding $u$:
        
        \begin{align}\label{eq:frequency_start}
            &-\rho \omega^2 \left(u - \frac13 e^2 \Delta u \right) - \frac{1}{2e}\operatorname{div}\operatorname{div}\mathcal{\hat{D}}\nabla\nabla u = \frac{1}{2e}f \\
            &\begin{cases}\label{eq:frequency_clamped_BC}
                u = g_\omega \\
                \frac{\partial u}{\partial n} = 0
            \end{cases} on\ \Gamma_c \\
            &\begin{cases}
                -\frac23 e^3 \rho\omega^2 \frac{\partial}{\partial n}u + \operatorname{div}\mathbf{M \cdot n} + \frac{\partial}{\partial \tau}(\mathbf{Mn \cdot \tau}) = 0 \\
                \mathbf{Mn \cdot \tau} = 0
            \end{cases} on\ \Gamma_f 
            \label{eq:frequency_fin}
        \end{align}
        
        For the practical computation, we suggest that the vibration stand, to which the specimen is clamped to, moves as a rigid body, thus the amplitude of boundary displacement $g$ does not depend on the spatial coordinate, $g(x,y) \equiv g$. 
        Let us denote the displacement amplitude in the test point, caused by the vibration with amplitude $g$, by $u_0$.
        Due to vibrations being harmonic, the amplitudes of the acceleration can be expressed as $-\omega^2 g$ at the boundary and as $-\omega^2 u_0$ at the test point.
        Additionally, due to the linearity of the problem, the acceleration caused by vibration with boundary displacement amplitude equal to unity, will be equal to
        \[
            u = \frac{u_0}{g} = \frac{-\omega^2 u_0}{-\omega^2 g}
        \]
        We suggest that the ratio of the complex amplitudes of the driving vibration and response $\frac{-\omega^2 u_0}{-\omega^2 g}$ can be found by analysing the experimental data.
        Thus, in the numerical solver, we fix $g \equiv 1$, solve the problem for the set of measured frequencies and compare it to the aforementioned ratio.
        We also suggest that the external loading, if present, is applied on the same frequency.
        
        \subsection{Formulation of the inverse problem}
        Let $(x_t, y_t)$ be the coordinates of the test point, i.e. the point above which the rangefinder is located.
        Let $\mathbf{\theta} = (\theta_1,\ \dots,\ \theta_k)$ be a set of model parameters, which uniquely define $D_{ij}$ and $\beta_{ij}$.
        For example, we can assume the isotropy of the sample by setting $\theta = (D, \nu, \beta)$ with 
            $\beta_{ij} = \beta$ and $D_{ij}$ defined according to the formula \eqref{eq:isotropicHooke}.
        Scaling and shifting of the variables can also be handled this way.

        The inverse problem is then formulated as follows. Given a set of driving frequencies $\{\omega_k,\ k \in \overline{1,\ N_\omega} \}$ 
            and the values of the complex amplitudes in the test point $\{u^{exp}_k,\ k \in \overline{1,\ N} \}$, acquired in the experiment,   
            find $\theta $ such that $u(x, y, \omega_k)$ is the solution of the BVP \eqref{eq:frequency_start}--\eqref{eq:frequency_fin} 
                with $\omega = \omega_k$  
            and $u^{exp}_k = u(x_t, y_t, \omega_k)$.

\section{Numerical method for direct problem}
    \subsection{Finite Element Solution}
        The weak form of the problem \eqref{eq:frequency_start}--\eqref{eq:frequency_fin} is\footnote{%
        TODO: for the article~---~accurately define the relevant functional spaces}
        \begin{equation}\label{eq:weak_form}
            \int\limits_{\Omega} \left(-\rho\omega^2(uv + 
                        \frac13 e^2 \nabla u \nabla v) + 
                        \frac1{2e}\mathcal{\hat{D}} \nabla\nabla u : \nabla\nabla v 
                        - fv \right)d\Omega = 0
        \end{equation}
        The solution is then approximated as a linear combination of basis functions $h_i$.
        In particular, Morley finite elements \cite{morley1968triangular} for the plate problem are used. 
        \begin{equation*}
            u = \sum\limits_{i \in I}u_ih_i + \sum\limits_{k \in D}g_k h_k
        \end{equation*}
        Here $I$ is a subset of indices of those basis functions, that are zero on $\Gamma_c$.
        $g_k = g(x_k, y_k)$ are set to approximate the Dirichlet boundary conditions \eqref{eq:frequency_clamped_BC} on $\Gamma_c$ and $u_i$ are the unknown coefficients which are to be found.
        Note that the boundary condition \eqref{eq:frequency_fin} is automatically satisfied by the solution.

        We arrive at the equation for $u_i$ by requiring that \eqref{eq:weak_form} holds for every $v \in \left\{h_i,\ i \in I \right\}$:
        \begin{equation}\label{eq:num_direct}
            \tilde{K}(\omega, \theta)u = \tilde{f}(\omega, \theta)
        \end{equation}
        with matrices in \eqref{eq:num_direct} defined as
        \begin{align}
            \tilde{K}(\omega, \theta) &= -\rho\omega^2\left(M + \frac13 e^2 L\right) + \sum\limits_{\alpha} D_{\alpha}(\theta)(1 + i\beta_{\alpha}(\theta))K^{\alpha}\\
            \tilde{f}(\omega, \theta) &= f_l -\rho\omega^2\left(f_M + \frac13 e^2 f_L \right) + \sum\limits_{\alpha} D_{\alpha}(\theta)(1 + i\beta_{\alpha}(\theta))f^{\alpha}
        \end{align}
        {\centering
        \begin{minipage}{0.49\textwidth}
            \begin{align}
                \left[M \right]_{ij} &= \int\limits_{\Omega} \rho h_ih_j d\Omega \\
                \left[L \right]_{ij} &= \int\limits_{\Omega} \rho (h_{i,x}h_{j,x} + h_{i,y}h_{j,y})d\Omega \\
                \left[K^\alpha \right]_{ij} &= \int\limits_{\Omega}V^\alpha(h_i, h_j)d\Omega 
            \end{align}
        \end{minipage}
        \hfill
        \begin{minipage}{0.49\textwidth}
            \begin{align}
                \left[f_l \right]_{j} &= \int\limits_{\Omega}f h_j d\Omega \\
                \left[f_M \right]_{j} &= -\sum\limits_{k \in D}g_k \int\limits_{\Omega} \rho h_k h_j d\Omega \\
                \left[f_L \right]_{j} &= -\sum\limits_{k \in D}g_k \int\limits_{\Omega} \rho (h_{k,x}h_{j,x} + h_{k,y}h_{j,y})d\Omega \\
                \left[f^\alpha \right]_{j} &= -\sum\limits_{k \in D}g_k\int\limits_{\Omega}V^\alpha(h_k, h_j)d\Omega 
            \end{align}
        \end{minipage}
        }
        and the bilinear forms $V^\alpha$ correspond to each elastic modulus and are defined as
        \begin{align}
            V^{11}(u, v) &= u_{xx}v_{xx} \\ 
            V^{12}(u, v) &= u_{yy}v_{xx} + u_{xx}v_{yy}\\ 
            V^{16}(u, v) &= u_{xy}v_{xx} + 2 u_{xx}v_{xy}\\ 
            V^{22}(u, v) &= u_{yy}v_{yy} \\ 
            V^{26}(u, v) &= u_{xy}v_{yy} + 2 u_{yy}v_{xy}\\ 
            V^{66}(u, v) &= 2 u_{xy}v_{xy}  
        \end{align}

        By the definition of the FE basis, some of the $u_j$ are equal to the value of $u$ at the nodes of the mesh, and some of them are equal to the values of $u_x$ or $u_y$.
        The test point, in general, may not be one of the nodal points, so interpolation is needed:\
        \begin{equation}
            u(x_t, y_t) = \sum\limits_{i \in I}h_i(x_t, y_t)u_i + \sum\limits_{k \in D}h_k(x_t, y_t)g_k = c^Tu + c_0 = P(u)
        \end{equation}

\begin{figure}
            \noindent\centering{
            \includegraphics[width=\textwidth]{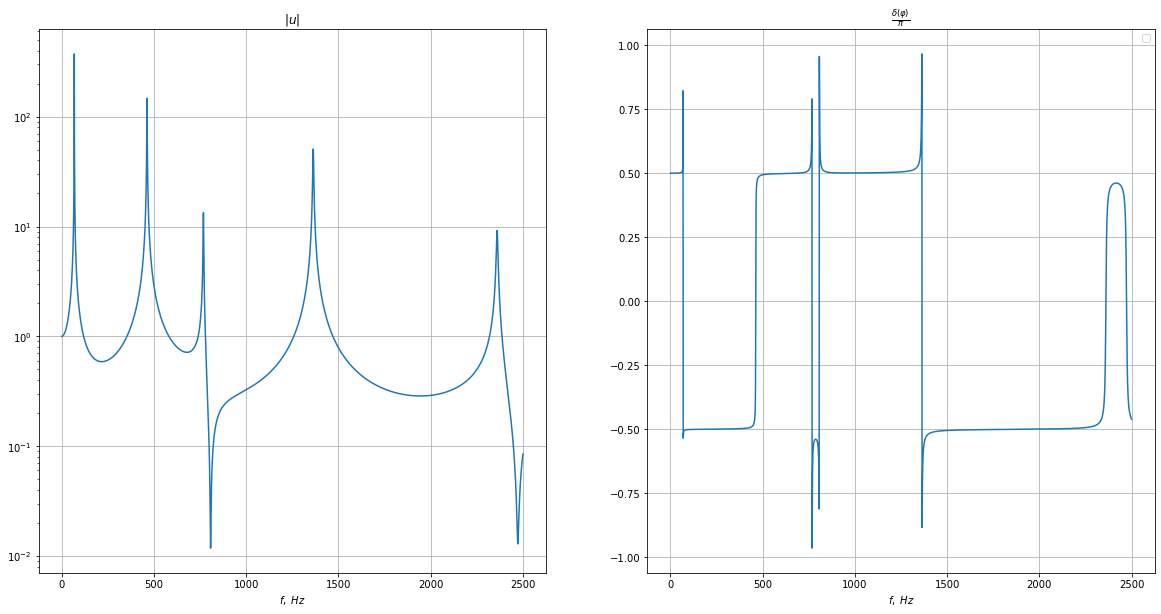}}
            \caption{Example AFC for the experimental setup, described below. Ratio of amplitudes and phase shift between driving vibration and response is plotted on the right and on the left, respectively.}\label{fig:sample_afc}
        \end{figure}

        By solving the problem \eqref{eq:num_direct} for different values of $\omega$ it is possible to numerically simulate the AFC of the sample, which can be seen on \figref{fig:sample_afc}.
        
        Note that the natural frequencies and decay factors, i.e. the values $\omega_k,\ \gamma_k$ such that the solution of the initial-value problem with zero boundary conditions is decomposed as a series of decaying vibrations $\sum\limits_{k = 1}^{+\infty} u_k e^{(-\gamma_k + i\omega_k)t}$ can also be approximately evaluated as the solution of the eigenvalue problem
        \begin{gather}\label{eq:natural_fs}
            \left(\sum\limits_{\alpha} D_{\alpha}(\theta)(1 + i\beta_{\alpha}(\theta))K^{\alpha} - \Lambda_k \cdot \left(M + \frac13 e^2 L\right)\right)u_k = 0\\
            \omega_k = \operatorname{Re}{\sqrt{\Lambda_k}},\ \gamma_k = \operatorname{Im}{\sqrt{\Lambda_k}}
        \end{gather}
        
        \subsection{Correction of inertia terms due to accelerometer mass}\label{sec:acc_mass_correction}
        If the specimen tested is rather small, the influence of the attached accelerometer cannot be neglected. 
        In our model, we suppose that the accelerometer mass is uniformly distributed through thickness of the specimen in the area, occupied by the accelerometer.
        Thus, the density $\rho$ in the weak form of the problem \eqref{eq:weak_form} is considered to be a function of coordinates, having the form
        \begin{equation}
            \rho(x, y) = \rho_0 + \chi_{\Omega_a}(x, y)\rho_{c}
        \end{equation}
        with $\rho_0$ being the specimen's material density, $\chi_{\Omega_a}$ is the characteristic function of the area $\Omega_a$, occupied by the accelerometer and $\rho_c = \frac{m_{a}}{2e S_a}$, $m_a, S_a$ being the mass and area of the accelerometer. 
        Consider the example of the mass matrix $M$:
        \begin{align}
            M &= \rho M_0 + \rho_c M_c \\
            \left[M_0 \right]_{ij} &= \int\limits_{\Omega} h_ih_j d\Omega \\
            \left[M_c \right]_{ij} &= \int\limits_{\Omega} \chi_{\Omega_a} h_ih_j d\Omega 
        \end{align}
        Thus, the procedure can be viewed as correction of a standard FE mass matrix $\rho M_0$ by a special matrix $\rho M_c$, defined above.

    \section{Optimization problem}
        We formulate the inverse problem as a nonlinear least-squares problem.
        Let $u(\theta, \omega)$ be the solution of \eqref{eq:num_direct}, $t$~---~the index of the test point.
        The process of finding the approximate solution to the inverse problem is formulated as minimization of the loss functional, which is defined as the mean-square error between the numerical solution an the experimental values.
        \begin{gather}
            \min_\theta L(\theta) = \frac{1}{N_\omega}\sum\limits_{k = 1}^{N_\omega}\left\| P(u(\theta, \omega_k)) - u^{exp}_k \right\|^2 \label{eq:loss} \\
                s.t.\ \theta > 0 \notag
        \end{gather}

        In the current work we suppose that the \emph{complex} amplitude, or, equivalently, both the amplitude ($|u^{exp}|$) and the phase shift ($\arg u^{exp}$) are known experimentally.
        Measurement of the phase shift between the accelerometers' signals requires good synchronization between the devices, which can be difficult, but the possibility of acquiring the parameters only from $|u^{exp}|$ is left for future studies.
        
        For the reconstruction of the parameters, we investigate two possible approaches:
        In the first approach, which will be referred to as local optimization, a <<sufficiently good>> initial guess for the parameters is known, and we want to improve this guess with fast, locally convergent methods.
        In the section below, we will try to define the requirements for the initial guess from numerical experiments.
        In the second approach, more realistic for the experimental composite material in question, only the orders of magnitude of each parameter can be estimated:
        \[
            \theta_i \in [\theta_i^L, \theta_i^U]; \quad \frac{\theta_i^U}{\theta_i^L} \sim 10^2-10^3
        \]
        This approach will be referred to as global optimization. 
        In this approach, we will use a heuristic global optimization method to obtain a set of parameters in the vicinity of the optimal solution in a reasonable amount of iterations, and then investigate the possibility of improving this result by means of local optimization.
    \subsection{Automatic differentiation}
        The simplicity of the underlying numerical model (the problem of evaluating $u(\theta, \omega)$ is reduced to the solution of a linear equation, the matrix and right-hand side are linear combinations of the unknown parameters with constant matrices) allows for usage of \emph{automatic differentiation methods} for fast and stable evaluation of the derivatives of the loss functional.
        In the present work, the \texttt{jax} library\cite{jax2018github} is used.
        The package provides tools for writing differentiable code that will automatically be just-in-time compiled for execution either on CPU or GPU.
        The theoretical studies show that, if $f(x)$ is an elementary scalar-valued function, the backward differentiation algorithm allows to compute the gradient of the function in no more than $3$ times than the time required for the computation of the function itself\cite{evtushenko_fad}.
        In our numerical experiments, run  with the usage of GPU, the evaluation of the loss function $L(\theta)$ and its derivatives up to the second order, takes around $1.5$ times longer than the function evaluation, which is really promising (although one should not forget that performance times can vary significantly between different machines, compilers, versions of libraries, etc.).
        A possibility for relatively fast evaluation of the derivatives, and the small dimensionality of the parameters space, is an argument for usage of second-order and quasi-newton first-order algorithms.
        These will be described in the following section.
    \subsection{Trust-region methods}
        It was chosen to implement the numerical methods for the optimization problem in the framework of \emph{trust-region methods}.
        The underlying idea is that during the minimization of the function $f(x)$, at each point $x_k$, the function is locally  substituted with its quadratic model
        \[
            m_k(p) = f_k + g_k^Tp + \frac12 p^TB_kp,\quad p = x - x_k,
        \]
        where $f_k = f(x_k),\ g_k = \nabla f(x_k)$ and the symmetric matrix $B_k$ is either the Hessian or some reasonable approximation of it. 
        The method also keeps track of a value $\Delta_k$, defining the radius of the ball, in which we believe the model to adequately represent the target function. 
        On each iteration, a constrained quadratic optimization problem is solved to find the direction to the new point:
        \begin{gather}\label{eq:tr_subproblem}
            p_k = \arg \min_{p} f_k + g_k^Tp + \frac12 p^TB_kp \\
            s.t.\ \|p_k\| \leq \Delta_k \notag
        \end{gather}
        After that, the improvement of the function, $f(x_k) - f(x_k + p_k)$ is compared with the improvement, predicted by the quadratic model.
        Depending on the relative improvement, the method decides if to accept the new point, and also adjusts the radius of the trust region. 
        The pseudo-code of the algorithm is presented below in the Algorithm \ref{alg:trust_region}.
        \begin{algorithm}[ht]
        \caption{Trust-region method}
        \label{alg:trust_region}
        \begin{algorithmic}[1]
            \REQUIRE $\mathbf{x}_0$~---~initial guess, $\Delta_{max} > 0,\ \Delta_0 \in (0, \Delta_{max}),\ \eta \in [0, \frac14)$
            \REPEAT
            \STATE Evaluate $p_k$ as the solution of \eqref{eq:tr_subproblem}
            \STATE $\rho_k = \frac{f(x_k) - f(x_k + p_k)}{m_k(0) - m_k(p_k)}$ \COMMENT{relative improvement}
            \IF{$\rho_k < \frac14$}
                \STATE $\Delta_{k+1} = \frac14\Delta_k$
            \ELSE
                \IF{$\rho_k > \frac{3}{4}$ \AND $ \|p_k\| = \Delta_k$}
                    \STATE $\Delta_{k+1} = \min(2\Delta_k,\ \Delta_{max})$
                \ELSE
                    \STATE $\Delta_{k+1} = \Delta_k$
                \ENDIF
            \ENDIF
            \IF{$\rho_k > \eta$}
                \STATE $x_{k+1} = x_k + p_k$
                \STATE Update the quadratic model
            \ELSE
                \STATE $x_{k+1} = x_k$
            \ENDIF
            \UNTIL{checkStopCondition() \OR $k \geq k_{max}$}
        \end{algorithmic}
        \end{algorithm}
        \pagebreak
        
        The type of the update formula used to evaluate the matrix $B_k$ in \eqref{eq:tr_subproblem} defines the certain algorithm in the framework of trust-region methods. 
        In the preliminary numerical experiments, two types of updates were considered.
        First one, called \emph{Newton-Gauss} method, uses the precise Hessian: $B_k = \nabla^2 f(x_k)$.
        The second one, called \emph{Broyden–Fletcher–Goldfarb–Shanno, or BFGS} method starts with exact Hessian ($B_0 = \nabla^2 f(x_0)$), then updates this matrix to  approximate the Hessian at each iteration by using the information from the gradients:
        \begin{align*}
            y &= g_{k+1} - g_k \\
            s &= x_{k+1} - x_k \\
            B_{k+1} &= B_k + \frac{yy^T}{s^Ty} - \frac{s^TB_k^TB_ks}{s^TB_ks}
        \end{align*}
        More information on the trust-region methods can be found in \cite{nocedal2006numerical}. 
        The  usage the BFGS update formula in the trust-region framework is reported in \cite{yuan_tr_bfgs}.
        In our numerical experiments, possibly due to effective parallelization in the code that was automatically compiled for execution on the GPU by  \texttt{jax}, the speed of evaluation of the exact Hessian is not more than $1.5$ times slower than the computation by the BFGS formula, and its usage provides significant decrease in the number of iterations, thus the numerical experiments below focus on Newton-Gauss method.

    \subsection{Global optimization}
        For exploring a larger subset of parameters in the approach of global optimization, a heuristic differential evolution method \cite{storn1997differential} was chosen, and its pseudocode is presented as Algorithm \ref{alg:de}.
        
        \begin{algorithm}[ht]
        \small
        \caption{Differential evolution}
        \label{alg:de}
        \begin{algorithmic}[1]
            \REQUIRE $CR\in(0;1)$~---~crossover rate; $0\leq F_{min} < F_{max} \leq 2$~---~ bounds for mutation factor; 
                $x_l \leq x^0_i \leq x_u,\ i \in \overline{1, NP}$~---~initial population
            \REPEAT
            \STATE Evaluate best member $b^k = \arg \min_i f(x^k_i)$ 
            \STATE Sample $F$ from $U\left([F_{min}, F_{max}) \right)$
            \FORALL{$i \in \overline{1, NP}$}
                \STATE Sample $i_1 \neq i_2$ 
                \STATE $b' = b^0 + F(x^k_{i_1} - x^k_{i_2})$
                \STATE $[v]_t = [b']_t$ with probability $CR$, else $[v]_t=[x^k_i]_t$ for $t\in \overline{1, \dim(x)}$
                \IF{$f(v) < f(x_i^k)$}
                    \STATE $x^{k + 1}_i = v$
                \ELSE
                    \STATE $x^{k + 1}_i = x^k_i$
                \ENDIF
            \ENDFOR
            \UNTIL $k\cdot NP \leq N_{fev, max}$ \OR $\operatorname{std}(x^k_i) \leq \varepsilon \operatorname{mean}(x^k_i)$
        \end{algorithmic}
        \end{algorithm}
        
        The method mimics biological processes of mutation and crossing-over.
        The set of search space points $x^k_i$ mimics a population, each parameter $\left[x^k_i \right]$ being some feature of a member <<organism>>.
        On every iteration, a candidate specimen is generated, by randomly combining the features of two other specimen (crossing-over), multiplied by a random factor (mutation). 
        If the value of the target functional is less for the new specimen, it is considered <<more fit>> and succeeds its predecessor in the population.    
        The method is easy to implement, has a small amount of hyperparameters, and for each of them it is quite clear how they affect the qualitative behavior of the method. Increasing the value of the crossover rate $CR$ or moving the boundaries of the interval of mutation factor closer to $2$ increases the magnitude of difference between the newly generated and existing population members, leading to broader area visited at the cost of convergence speed. In our numerical experiments, implementation of DE provided by the \texttt{scipy.optimize} package is used. 

\section{Numerical experiments}
    \subsection{Problem statement for an isotropic strip}
        We start investigating the behaviour of our method from a simple test case, which is easy to implement in an experimental environment, and with material parameters known from literature.
        The experimental specimen is a steel strip, clamped on one of the short sides, with the accelerometer being attached in the proximity of the second short side.
        The physical parameters are compiled in the table~\ref{tab:strip_params}.
        \begin{table}[ht]
            \centering
            \begin{tabular}{l|c}
                Dimensions, $[mm]$ & $100 \times 20 \times 1$ \\
                Density, $\left[\frac{kg}{m^3}\right]$ & $7920$ \\
                Young's modulus, $[GPa]$ & $198$ \\
                Shear modulus, $[GPa]$ & $77$ \\
                Poisson's ratio, $[1]$ & $0.286$ \\
                Flexural rigidity, $\left[Pa \cdot m^3\right]$ & $17.97$ \\
                Loss factor, $ [1]$ & $0.003$ \\
                Accelerometer mass, $g$ & $1$\\
                Accelerometer radius, $mm$ & $1$
            \end{tabular}
            \caption{Physical parameters of the isotropic test specimen}
            \label{tab:strip_params}
        \end{table}
        
        \begin{figure}
            \centering
            \begin{subfigure}[ht]{\textwidth}
                \noindent\centering{
                \includegraphics[width=\textwidth]{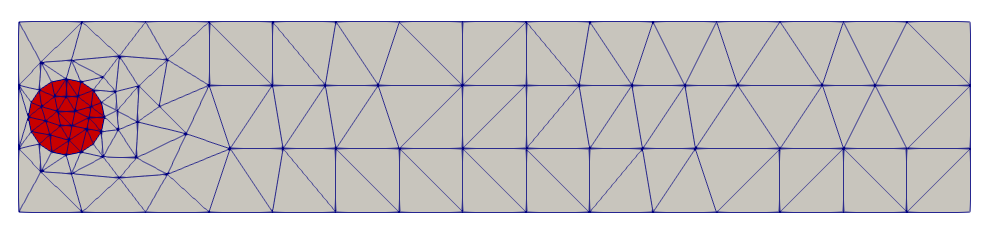}}
                \subcaption{Symmetrical placement of the accelerometer}
                \label{fig:mesh_sym}
            \end{subfigure}
            \begin{subfigure}[hb]{\textwidth}
                \noindent\centering{
                \includegraphics[width=\textwidth]{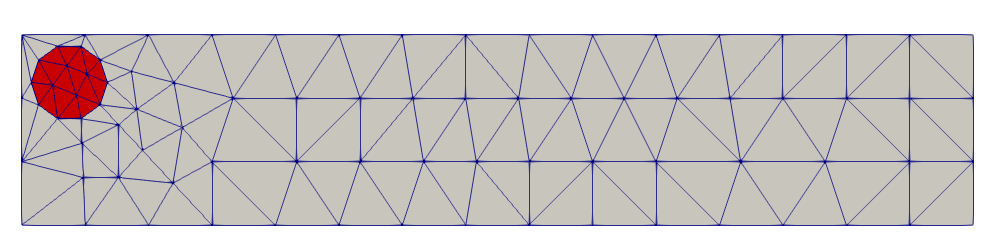}}
                \subcaption{Accelerometer shifted towards the longer side. This setup was chosen for final experiments}
                \label{fig:mesh_shift}
            \end{subfigure}
            \caption{Mesh for the numerical experiments with the isotropic strip. The right side is clamped, and the area, occupied by the accelerometer is highlighted in red.}
        \end{figure}
        In the preliminary direct problem computations, we would like to examine the influence of the method of accounting for the accelerometer's mass, suggested in Section~\ref{sec:acc_mass_correction}, and of different positioning of the accelerometer.
        \begin{figure}[hb]
            \noindent\centering{
            \includegraphics[width=\textwidth]{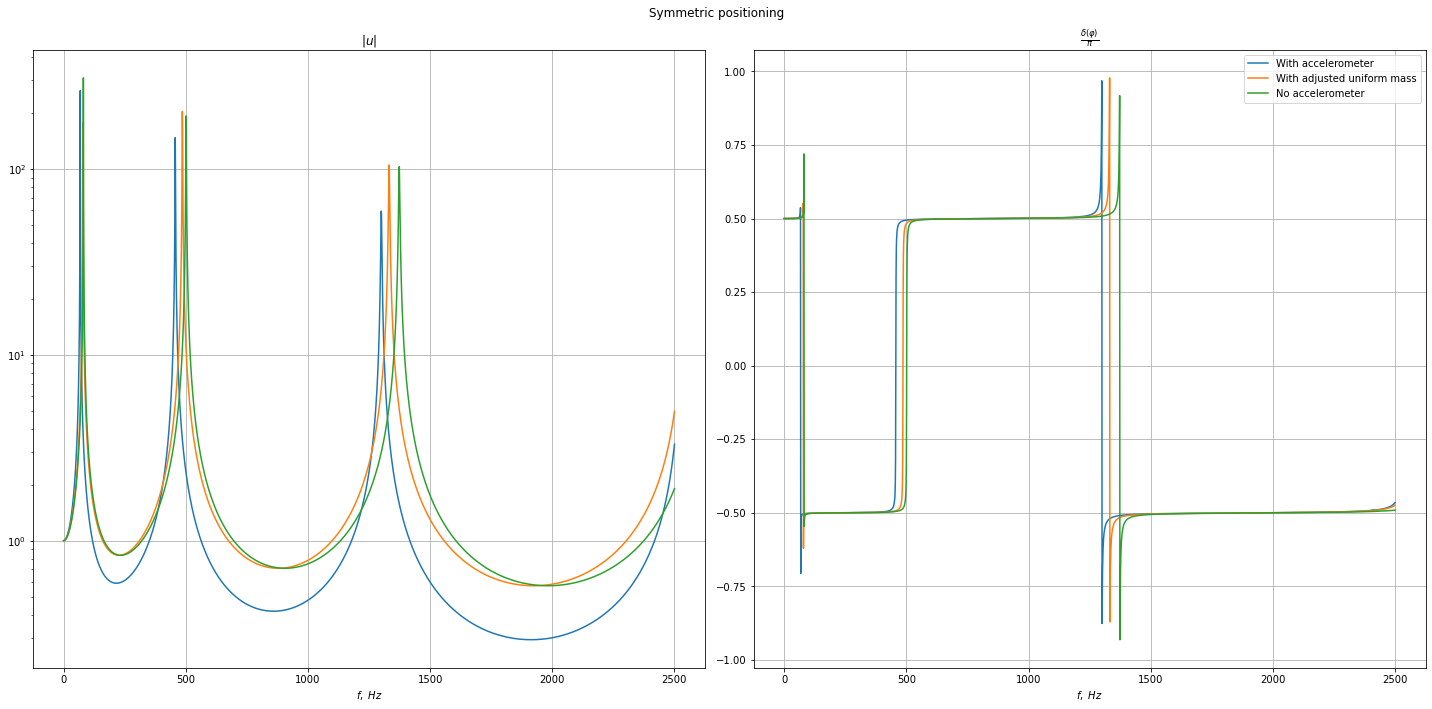}}
            \caption{AFCs, obtained with different approaches for modelling the influence of the accelerometer.}
            \label{fig:accel_corrections}
        \end{figure}
        The AFCs, computed for the case with symmetric positioning of the accelerometer (\figref{fig:mesh_sym}), are presented of \figref{fig:accel_corrections}.
        The AFCs computed with correction for the influence of the accelerometer is plotted in blue. 
        The approach is compared with two alternative ones. 
        First one, in orange is evaluated as if there was no accelerometer at all.
        For the second one, in green, density of the material is increased in such a way that the effective mass is equal to the total mass of the specimen and the accelerometer. 
        It is clearly visible, that the positioning of peaks and their width significantly changes due to the introduced correction, and neglecting the accelerometer's mass and position can provide incorrect results.
        We thus adhere to using this method in the following computations.
        
         \begin{figure}[ht]
            \noindent\centering{
            \includegraphics[width=\textwidth]{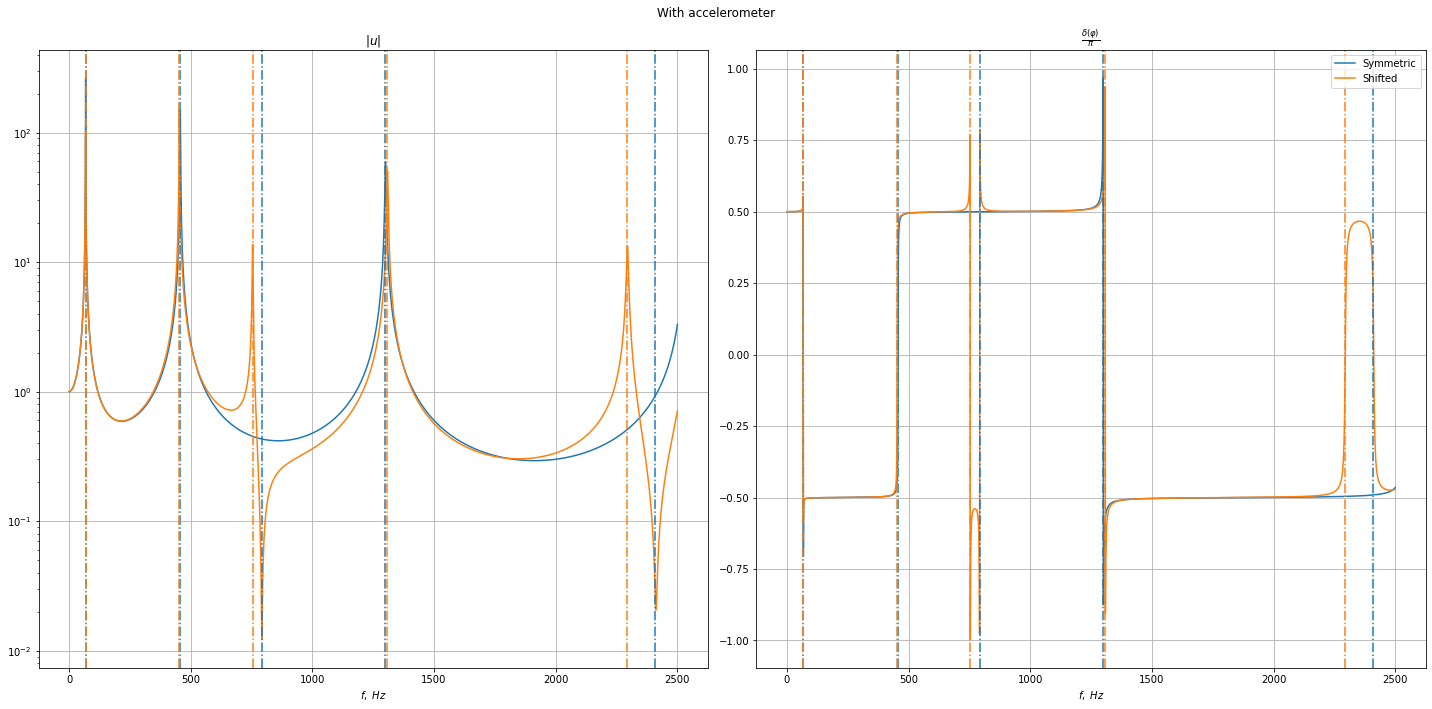}}
            \caption{AFCs for symmetric and shifted positioning of the accelerometer. }
            \label{fig:accel_sym_v_shift}
        \end{figure}
        We would also like to examine the effect of non-symmetric positioning of the accelerometer.
        The mesh for the test case, with the axis of the accelerometer shifted $5~mm$ from the center line towards the short side, is presented on \figref{fig:mesh_shift}.
        The AFCs for the symmetric and shifted variants of positioning are presented on \figref{fig:accel_sym_v_shift}, in blue and orange, respectively.
        We also represent natural frequencies, defined by \eqref{eq:natural_fs}, with dash-dotted vertical lines of the same colors.
        The first, second and fourth of the natural frequencies are quite close to each other, and correspond to an extremum in the respective AFC.
        On the other hand, there the peak in the proximity of the third and the fifth natural frequency appear only for the case for the case with non-symmetric positioning.
        Our preliminary experiments in global optimization show that larger amount of such peaks in the examined frequency range improve the chances of finding the global minimum.
        Thus, in the following computational experiments the test case with shifted accelerometer \figref{fig:mesh_shift} is studied.
        
        In the following sections, four sets of frequencies are considered, with $201$ equidistant points ranging from zero to $f_{max} \in [200, 600, 1000, 1500]~Hz $.
        There are $1,\ 2,\ 3$ and $4$ peaks on the ranges respectively.
        Reference data is generated by solving the forward problem with reference parameters from Table~\ref{tab:strip_params}, and, possibly, added noise, which was sampled from normal distribution with zero mean and variance, proportional to the maximal absolute value in the AFC. 
        The proportionality constant is referred to as noise level and expressed in percent.
        The quality of the final solution is assessed based on the relative error of each parameter, defined by the formula
        \begin{equation}
            re(\theta_i) = \frac{\theta_i - \theta_i^{ref}}{|\theta_i^{ref}|}
        \end{equation}

    \subsection{Local optimization example}
        The parameters we use for local optimization are 
        \begin{equation}\label{eq:param_transform_isotropic}
            \theta_1 = D,\ \theta_2 = \nu,\ \theta_3 = \beta
        \end{equation}
        Here flexural rigidity $D$ and Poisson's ratio $\nu$ define the compliance matrix as in \eqref{eq:isotropicHooke}, and, due to isotropy of the material, all the loss factors are equal: $\beta_\alpha = \beta$.
        It can be seen that the reference parameters have different orders of magnitude.
        Although it is reported in literature that the usage of scaling of parameters may improve the convergence of the method, in this particular case, the usage of shifted and scaled parameters in that way that the new parameters are of the same order of magnitude, didn't give any additional benefits.
        
\begin{figure}[t]
            \noindent\centering{
            \includegraphics[width=\textwidth]{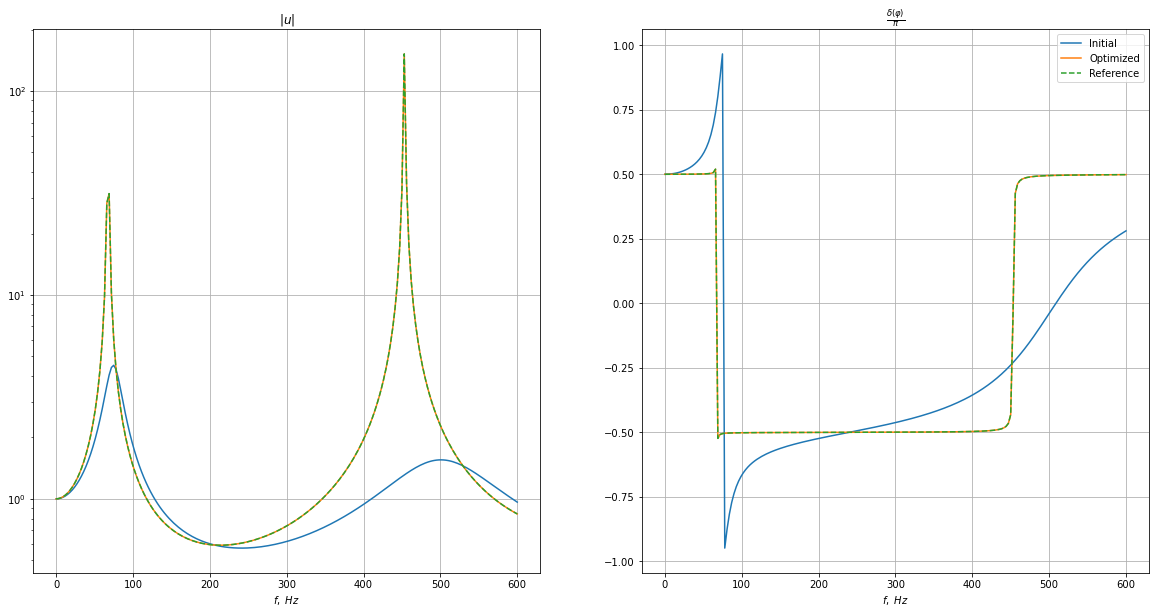}}
            \caption{AFCs before and after local optimization for $f_{max} = 600\ Hz$}
            \label{fig:local_600}
        \end{figure}

\begin{figure}
            \centering
            \begin{subfigure}[t]{.47\textwidth}
                \noindent\centering{
                \includegraphics[width=\textwidth]{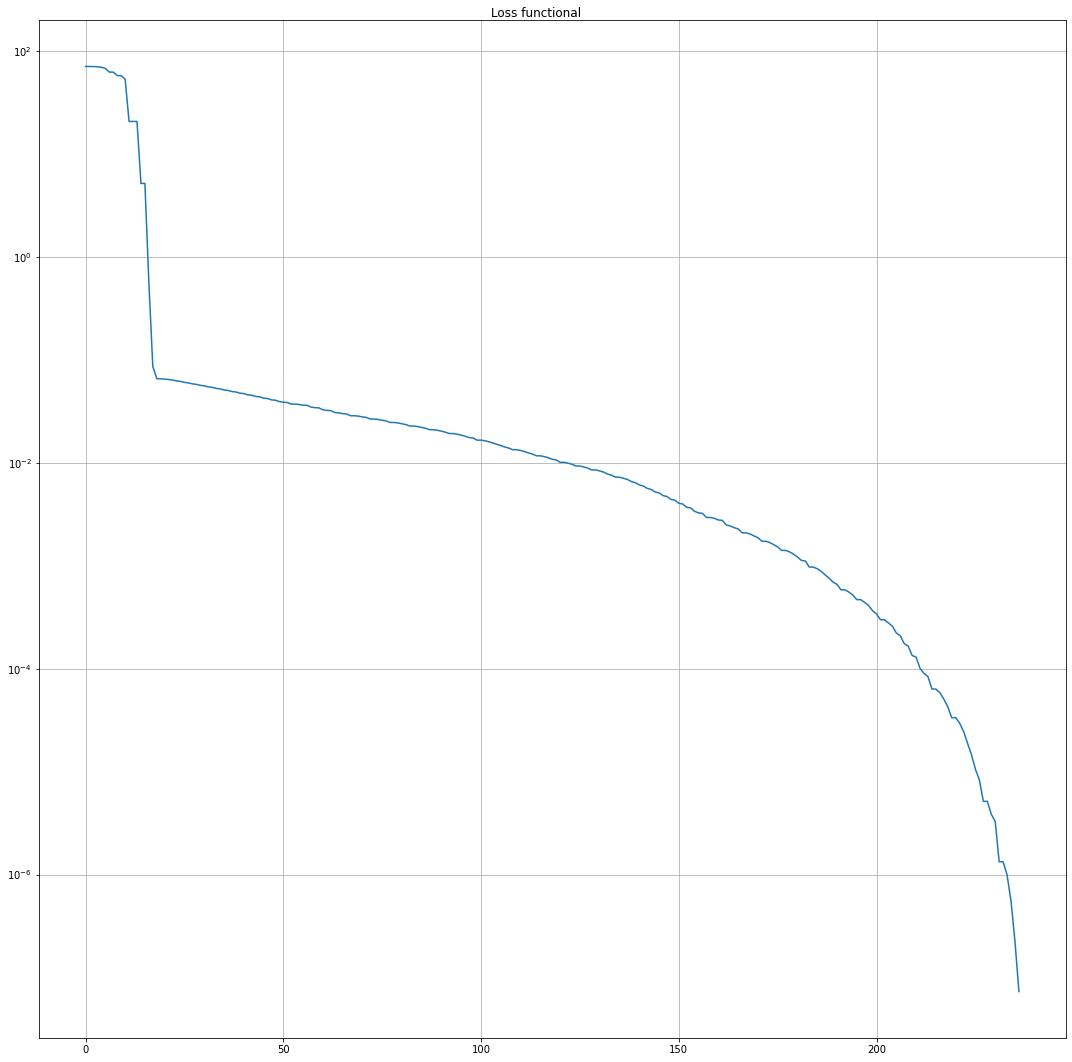}}
                \subcaption{Loss functional}
                \label{fig:local_600_loss}
            \end{subfigure}
            \hfill
            \begin{subfigure}[t]{.47\textwidth}
                \noindent\centering{
                \includegraphics[width=\textwidth]{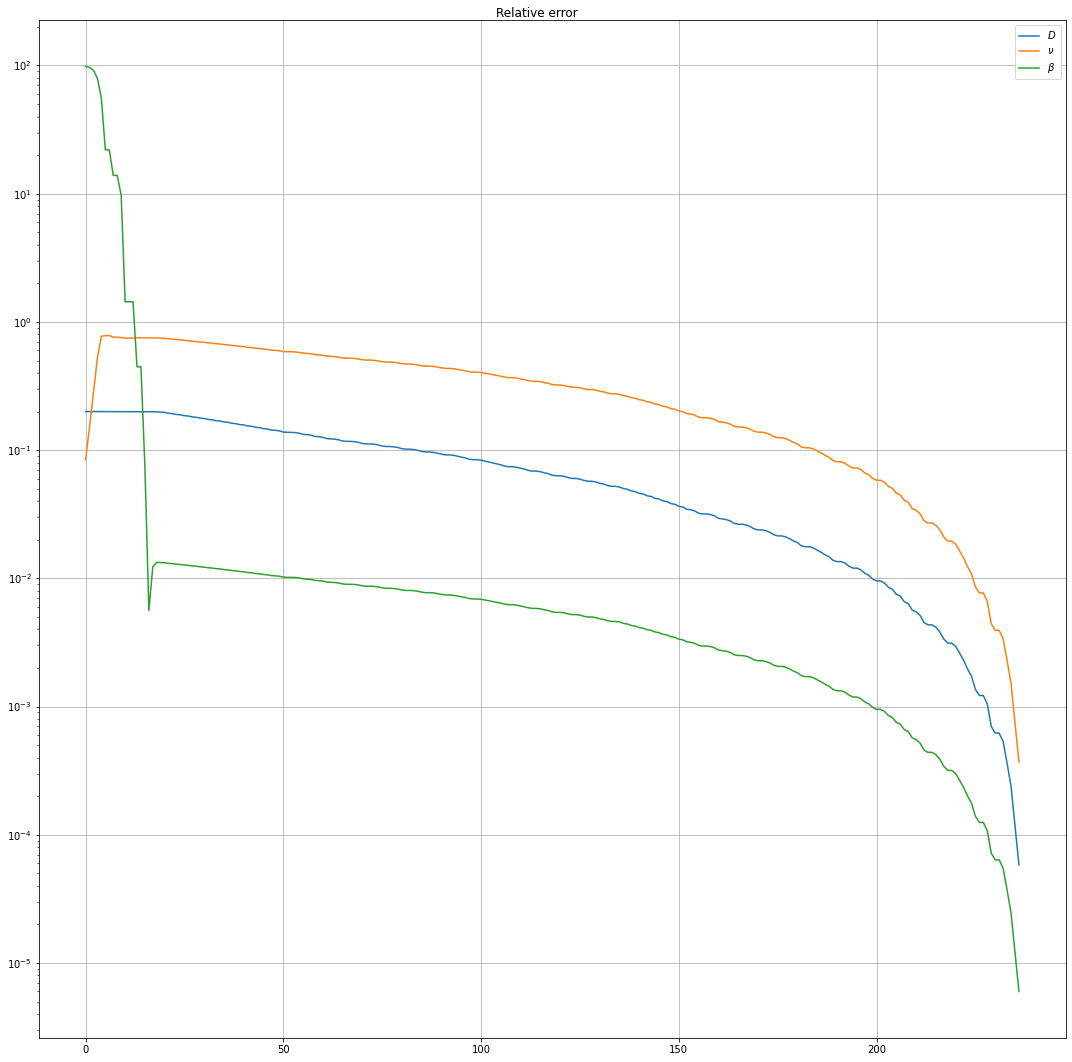}}
                \subcaption{Relative error}
                \label{fig:local_600_re}
            \end{subfigure}
            \caption{Convergence of the trust region method with initial relative error $(0.20,\ 0.05,\ 99.0)$}
            \label{fig:local_600_conv}
        \end{figure}

        Out of the three parameters, the loss factor $\beta$ can be reconstructed from the initial guesses. which significantly, up to two orders of magnitude, differ from the reference value.
        For the first two parameters, it is only possible to improve a rather precise initial guess.
        For example, in case of the frequency range $(0, 600)$, with a rather good initial guess in $\nu$ (relative error $5\%$), it was possible to solve the problem with initial error in $D$ of $20\%$ and initial $\beta$ one hundred times larger than the reference. 
        The final relative error was $(6\cdot 10^{-5},\ 4 \cdot 10^{-4},\ 6 \cdot 10^{-6})$.
        The AFCs for initial and final parameters, compared with the reference, can be seen on \figref{fig:local_600}, and the behaviour of the algorithm's loss functional and relative error is depicted on \figref{fig:local_600_conv}. Note that the trust-region method converges monotonically and superlinearily.  
        
\begin{figure}[t]
            \noindent\centering{
            \includegraphics[width=\textwidth]{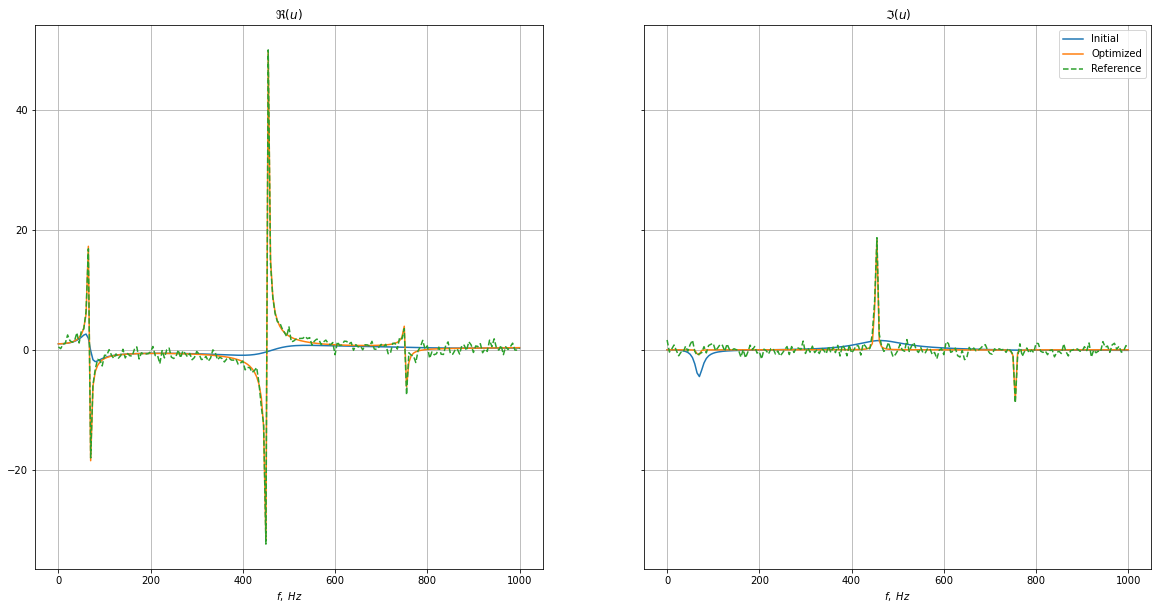}}
            \caption{AFCs before and after local optimization. Note that the real and imaginary part of the AFC are plotted instead of amplitude and phase shift. The green curve represents reference data with added Gaussian noise.}
            \label{fig:local_1000}
        \end{figure}

\begin{figure}
            \centering
            \begin{subfigure}[t]{.47\textwidth}
                \noindent\centering{
                \includegraphics[width=\textwidth]{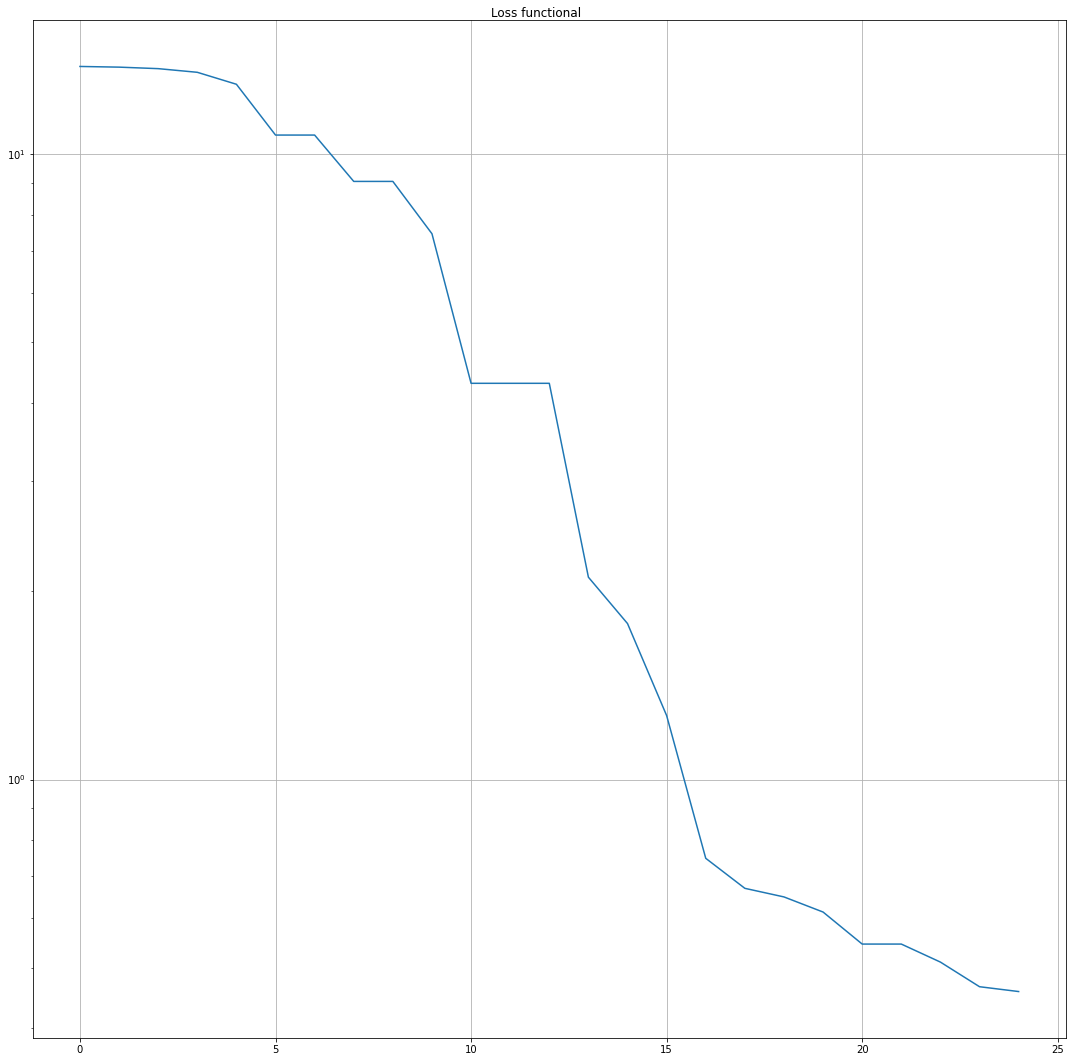}}
                \subcaption{Loss functional}
                \label{fig:local_1000_loss}
            \end{subfigure}
            \hfill
            \begin{subfigure}[t]{.47\textwidth}
                \noindent\centering{
                \includegraphics[width=\textwidth]{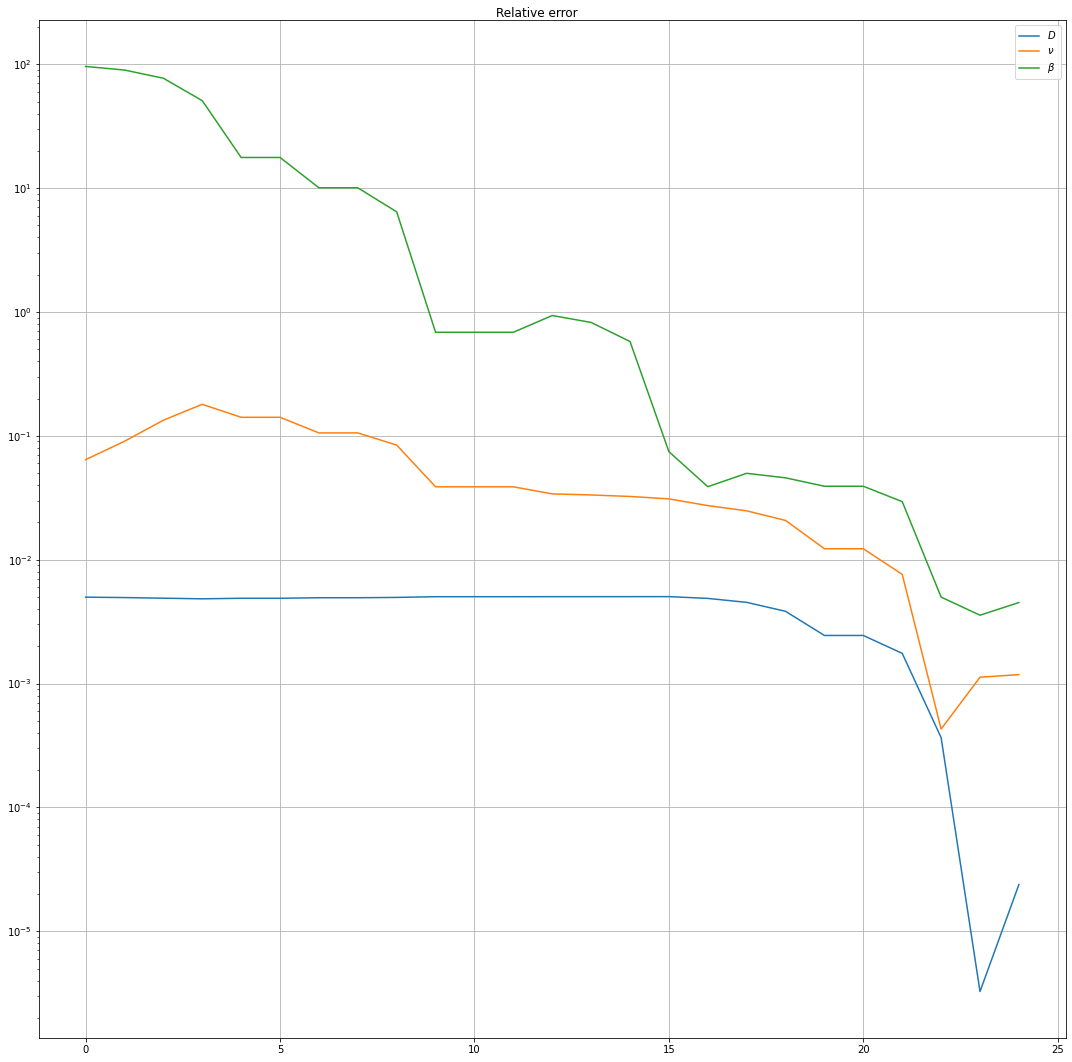}}
                \subcaption{Relative error}
                \label{fig:local_1000_re}
            \end{subfigure}
            \caption{Convergence of the trust region method with initial relative error $(0.005,\ 0.05,\ 99.0)$}
            \label{fig:local_1000_conv}
        \end{figure}

        In the second example, with AFCs and convergence results presented on Figures \ref{fig:local_1000} and \ref{fig:local_1000_conv} respectively, the local method is applied to a frequency range with three peaks, maximum frequency of $1000\ Hz$ and a moderate level of noise $1\%$. 
        For this case, a really good initial guess in $D$, with only $0.5\%$ difference from the reference value, was required. 
        It was then possible, with usage of the trust-region method, to improve the solution by orders of magnitude.
        The final relative error was $(2\cdot 10^{-5},\ 1 \cdot 10^{-3},\ 5 \cdot 10^{-3})$.
        Although the requirement on initial precision of $D$ might sound constraining, in the following, we show that generation of such initial guesses is, in principle, possible with usage of the global optimization method.
        
        \clearpage
        
    \subsection{Results of global optimization}
        Taking into account the properties of the local optimization process, noted in the section above, the following strategy for global minimization was chosen.
        As the loss factor $\beta$ can be stably improved by the local method, its value is fixed, equal to $0.01$, which is larger than the reference value $0.003$.
        As the order of magnitude of $D$ and $\nu$ are significantly different, we also scale this parameter, so they contribute equally to the stopping criterion of the algorithm.
        Thus, the parameters and their bounds are 
        \begin{equation}
            \theta_1 = D \in [1,\ 100],\ \theta_2 = \nu \cdot 100 \in [0,\ 0.5\cdot 100]
        \end{equation}
        The bounds on $\nu$ are physical restrictions on the Poisson's ratio and we believe that $D$ is known very vaguely, up to two orders of magnitude.
        The reference AFC is generated for the frequency range, mentioned above.
        For each of them, the reconstruction is studied for data without noise, as well as for data with $1\%$ and $3\%$ noise.
        
        The population size was set equal to $15 \cdot 2$.
        As for the parameters of the Differential Evolution method, the bounds for the mutation factor were chosen equal to $(0.7, 1.0)$, to increase the search space.
        The possible decrease in convergence speed is not critical, as we only generate a <<good enough>> initial guess to be improved by the local method.
        For the same reason we set $\varepsilon$ in the stopping criterion to a rather non-restrictive value of $10^{-2}$.
        Due to the stochastic nature of the method, it is sensitive to the initial population, thus, for every reference AFC, five restarts of the method with different initial population were done.
        Then, the result of the restart with minimal value of the loss functional is taken as the start parameters for the trust-region method.
        Note that for the local method, we use the parameter transformation, described in \eqref{eq:param_transform_isotropic}.
        
        \begin{table}[hb]
\centering
\begin{tabular}{ll|ccc|ccc}
\toprule
     &      &                 & R.E. after global&  &          &     R.E. after local &  \\
Noise   & $f_{max}$&   $D$ & $\nu$ & $\beta$ & $D$ & $\nu$ & $\beta$                                \\
\midrule
$0\% $  & $200 $  &    $1.3 \cdot 10^{-1}$ & $ 5.8 \cdot 10^{-1}$ & $ 2.3$ &    $1.1 \cdot 10^{-1}$ & $ 5.0 \cdot 10^{-1}$ & $ 5.6 \cdot 10^{-5}$ \\
        & $600 $  &    $1.8 \cdot 10^{-1}$ & $ 7.0 \cdot 10^{-1}$ & $ 2.3$ &    $1.8 \cdot 10^{-5}$ & $ 1.1 \cdot 10^{-4}$ & $ 1.8 \cdot 10^{-6}$ \\
        & $1000$  &   $-9.9 \cdot 10^{-4}$ & $ 2.5 \cdot 10^{-3}$ & $ 2.3$ &  $-1.3 \cdot 10^{-12}$ & $ -6.9 \cdot 10^{-12}$ & $ 8.2 \cdot 10^{-11}$ \\
        & $1500$  &    $3.8 \cdot 10^{-3}$ & $ 1.3 \cdot 10^{-2}$ & $ 2.3$ &   $2.5 \cdot 10^{-9}$ & $ 1.7 \cdot 10^{-8}$ & $ -6.4 \cdot 10^{-9}$ \\
\midrule
$1\% $  & $200 $  &    $1.9 \cdot 10^{-2}$ & $ 1.2 \cdot 10^{-1}$ & $ 2.3$ &    $2.0 \cdot 10^{-1}$ & $ 7.8 \cdot 10^{-1}$ & $ 8.0 \cdot 10^{-3}$ \\
        & $600 $  &    $2.0 \cdot 10^{-1}$ & $ 7.4 \cdot 10^{-1}$ & $ 2.3$ &  $-4.4 \cdot 10^{-2}$ & $ -3.5 \cdot 10^{-1}$ & $ 8.9 \cdot 10^{-3}$ \\
        & $1000$  &   $-1.4 \cdot 10^{-3}$ & $ 2.6 \cdot 10^{-4}$ & $ 2.3$ &  $-2.8 \cdot 10^{-5}$ & $ -1.1 \cdot 10^{-3}$ & $ 4.5 \cdot 10^{-3}$ \\
        & $1500$  &    $6.3 \cdot 10^{-3}$ & $ 3.0 \cdot 10^{-2}$ & $ 2.3$ &    $4.0 \cdot 10^{-3}$ & $ 2.6 \cdot 10^{-2}$ & $ 2.2 \cdot 10^{-2}$ \\
\midrule
$3\% $  & $200 $  &    $1.3 \cdot 10^{-1}$ & $ 5.8 \cdot 10^{-1}$ & $ 2.3$ &    $2.8 \cdot 10^{-1}$ & $ 9.4 \cdot 10^{-1}$ & $ 2.1 \cdot 10^{-2}$ \\
        & $600 $  &    $1.1 \cdot 10^{-1}$ & $ 5.1 \cdot 10^{-1}$ & $ 2.3$ &  $-6.6 \cdot 10^{-2}$ & $ 1.4$                & $ 2.1 \cdot 10^{-2}$ \\
        & $1000$  &  $-1.9 \cdot 10^{-3}$ & $ -2.2 \cdot 10^{-3}$ & $ 2.3$ &  $1.7 \cdot 10^{-4}$ & $ -1.7 \cdot 10^{-3}$ & $ -6.0 \cdot 10^{-3}$ \\
        & $1500$  &    $4.0 \cdot 10^{-3}$ & $ 1.7 \cdot 10^{-2}$ & $ 2.3$ &    $4.5 \cdot 10^{-3}$ & $ 3.0 \cdot 10^{-2}$ & $ 6.4 \cdot 10^{-2}$ \\
\bottomrule
\end{tabular}
\caption{Results of global optimization and polishing of local optimization}
\label{tab:opt_isotropic}
\end{table}

        The results of the global optimization process are presented in Table~\ref{tab:opt_isotropic}.
        As was suggested before the experiment, the presence of several peaks in the reference AFC is crucial for locating the global minimum.
        The reconstruction with reference data possessing only one peak ($f_{max} = 200\ Hz$) didn't lead to successful results, 
            and as for the frequency range with two peaks ($f_{max} = 600\ Hz$), a decent value of $\nu$ was only achieved on a case without noise.
        On the other hand, for the cases with three and four peaks ($f_{max} = 1000$ and $1500\ Hz$ respectively), the proximity of the global optimum can reliably be reached in at least one of five restarts, and for the cases without noise and moderate $1\%$ noise, the global result is significantly improved by usage of the trust-region method.
        It is worth noting that better results were achieved for the cases with $f_{max} = 1000$.
        As the same number of frequency points were used in all of the cases, we suggest that the smaller distance between the points may be beneficial for the local convergence of the algorithm.

    
    
\section{Conclusion}
    In the present work, the reasons for the development of a nondestructive method of acquring the elastic parameters of a thin composite plate from the data of the vibrational testing are explained.
    The experimental stand is described and the equations, governing the dynamics of the specimen during the testing are studied.
    The problem of acquiring the AFC of the specimen is formalised as a series of BVPs.
    A numerical method for solving these BVPs is provided.
    
    The inverse problem is formulated as a nonlinear least-square optimization problem, for which fast and stable evaluation of derivatives is possible with automatic differentiation methods. 
    The usage of trust-region Newton method in this case allowed for fast improvement of a decent initial guess, and numerical experiments with usage of the heuristic Differential Evolution have shown that generation of such initial guesses is possible with quite vague prior knowledge on the parameters.
    An important observation is that the presence of several peaks in the AFC on the studied frequency range was crucial for arriving at the correct solution.
    
    The suggested approach of using  differentiable programming in the solution of a direct problem can be extended for more complex tasks, such as accounting for frequency dependence of the material parameters, or using appropriate Fourier or Laplace transforms for solving the problems with time-domain data.
    As for the solution of the inverse problem, work can be done in several directions. 
    The efficiency of solution of linear systems on the GPU should be studied more carefully, and alternative iterative methods with joint optimization of the parameters and solution of the direct problem can be suggested. 
    The conditions for reliable global solution of the problem have to be studied more carefully, and design of better loss functionals than the standard mean-square error may be possible.

\section{Conclusion}
    The authors would like to thank Andrey Karchevsky and Maxim Shishlenin for meaningful feedback on the intermediate steps of the project.

\bibliographystyle{ieeetr}
\bibliography{bibliography}

\end{document}